\documentclass[11pt]{article}
\usepackage{amssymb}
\usepackage{amsfonts}
\usepackage{amsmath}
\usepackage{geometry}
\usepackage[colorlinks=true, linkcolor=red, citecolor=blue]{hyperref}
\usepackage{float}
\usepackage{mathtools}
\usepackage{caption}
\usepackage{graphicx}
\usepackage[font=small,labelfont=bf,labelsep=space]{caption}
\usepackage{subfig}
\usepackage{setspace}

\makeatletter

\@addtoreset{equation}{section}
\makeatother

\input{tcilatex}
\begin{document}

\author{Zouha\"{\i}r Mouayn \\
{\footnotesize \ Department of Mathematics, Faculty of Sciences and Technics
(M'Ghila),\vspace*{-0.2em}}\\
{\footnotesize \ Sultan Moulay Slimane University, BP. 523, B\'{e}ni Mellal
23000, Morocco }\\
{\footnotesize Institut des Hautes \'{E}tudes Scientifiques, Paris-Saclay
University,}\\
{\footnotesize 35 route de Chartres, 91893 Bures-sur-Yvette, France.}\\
{\footnotesize e-mail: mouayn@gmail.com }\vspace*{-0.5em}}
\title{Integrated density of states for the sub-Laplacian on Heisenberg
groups }
\maketitle

\begin{abstract}
We rederive the expression of the integrated density of states for \ the
sub-Laplacian on Heisenberg groups $\mathbb{H}_{n}$ by using its resolvent
kernel.
\end{abstract}

\section{Introduction}

The integrated density of states (IDS) has its origin in physics, that when
studying properties of the spectrum of a Hamiltonian operator $H$, it turns
out that one of the relevant objects is a distribution function, associated
to a certain probability measure, encoding the distribution of the spectrum
of $H$. Note that spectrum of a self-adjoint operator is a subset of the
real $\mathbb{R}.$ Hence, this distribution function is defined on $\mathbb{R%
}$ as well$.$ \ The IDS is the function\ defined by 
\begin{equation}\label{eq1}
\mathbb{R}\ni \lambda \rightarrow \ \mathcal{N}\left( \lambda \right)
:=\lim_{\Omega \rightarrow \mathbb{R}^{d},\text{ }\Omega \in \mathcal{O}}%
\frac{\sharp \left\{ \lambda _{j}\leq \lambda \right\} }{Vol(\Omega )} 
\tag{1.1}
\end{equation}%
where $\sharp \left\{ \lambda _{j}\leq \lambda \right\} $\ is the number of
eigenvalues less than $\lambda $\ of $H$\ restricted to a bounded domain $%
\Omega $\ with appropriate boundary conditions. $\mathcal{O}$\ is a family
of bounded open sets which is introduced to specify the way expanding $%
\Omega $\ to whole space\ $\mathbb{R}^{d}.$ From the definition \eqref{eq1}, we see that the integrated density of states gives the average
number of the possible energy levels below the energy $\lambda $\ per unit
volume.

\bigskip

One may also define the IDS by (\cite{AS},\cite{HS}): 
\begin{equation}
\mathcal{N}\left( \lambda \right) :=\lim_{\Omega \rightarrow \mathbb{R}%
^{d},\Omega \in \mathcal{O}}\frac{tr(\chi _{\Omega }P_{\left( -\infty
,\lambda \right] }\left( H\right) \chi _{\Omega })}{Vol(\Omega )}  \tag{1.2}
\end{equation}%
where $\chi _{\Omega }$\ denotes the multiplication operator by the
characteristic function of the domain $\Omega \subset \mathbb{R}^{d}$ and $%
P_{J}\left( H\right) $\ denotes the spectral projection of $H$\ for a Borel
set $J\subset \mathbb{R}.$ In this sense the IDS is an average of the kernel
of the spectral projection operator and $\cite{Stri}$ have showed that it
can serve as a substitude for the eigenvalue counting function $\sharp
\left\{ \lambda _{j}\leq \lambda \right\} $. This definition is motivated by
well known Weyl asymptotic law%
\begin{equation}\label{eq13}
\sharp \left\{ \lambda _{j}\leq \lambda \right\} =c_{n}Vol\left( \Omega
\right) \lambda ^{\frac{1}{2}n}+o\left( \lambda ^{\frac{1}{2}n}\right)  
\tag{1.3}
\end{equation}%
when $H=-\Delta $, the Laplacian on a compact domain $\Omega \subset \mathbb{%
R}^{d}$ with smooth boundary with Dirichlet or Neumann boundary conditions.
\ By dividing \eqref{eq13} by the volume of $\Omega $ and observing
that the ratio $N\left( \lambda \right) /Vol\left( \Omega \right) $ is equal
to the average value of the diagonal function $\Phi \left( \lambda
;x,x\right) $ on $\Omega ,$ where $\Phi \left( \lambda ;x,y\right) $ stands
for the kernel of the spectral projection operator $E_{\lambda }$ onto the $%
\left[ 0,\lambda \right] $ of the spectrum,%
\begin{equation}
\Phi \left( \lambda ;x,y\right) =\sum\limits_{\lambda _{j}\leq \lambda }\phi
_{j}\left( x\right) \overline{\phi _{j}\left( y\right) }  \tag{1.4}
\end{equation}%
where $\left\{ \phi _{j}\right\} $ is an orthonormal basis of eigenfunctions
associated with eigenvalues $\left\{ \lambda _{j}\right\} .$ Indeed, we have 
\begin{equation}\label{eq15}
\mathcal{N}\left( \lambda \right) =\frac{1}{Vol\left( \Omega \right) }%
\int\limits_{\Omega }\Phi \left( \lambda ;x,x\right) d\mu \left( x\right) , 
\tag{1.5}
\end{equation}%
where $d\mu $ denotes the Lebesgue measure or the Riemannian measure on $%
\Omega $ and one may rewrite \eqref{eq13} as follows%
\begin{equation}
\sharp \left\{ \lambda _{j}\leq \lambda \right\} =c_{n}\lambda ^{\frac{1}{2}%
n}+o\left( \lambda ^{\frac{1}{2}n}\right) ,  \tag{1.6}
\end{equation}%
$c_{n}$ is constant that depends on the dimension and can be explicitly
determined.

\medskip

Since $\mathcal{N}\left( \lambda \right) $ is nondecreasing function, its
derivative $\rho \left( \lambda \right) =d\mathcal{N}\left( \lambda \right) $%
\ $/d\lambda $\ is well defined as Steiljes-Lebesgue positive measure and is
called the density of states (DOS). The latter can be exploited, in the case
of the Landau Hamiltonian, to measure, in the gap between Landau levels,
energies of new states that could be created if a random potential is added $%
\cite{BVS}.$ For different mathematical interest on the IDS we refer to $%
\cite{KM}$ and references therein.

\smallskip \smallskip 

In this Note, our aim is to rederive the expression of the integrated
density of states for the sub-Laplacian on Heisenberg groups $\mathbb{H}_{n}$
by using its resolvent kernel as obtained in \cite{AskMou}.

\bigskip

\section{Magnetic Laplacians in $\mathbb{C}^{n}$}

In suitable units and up to additive constant, the motion of a charged
particle in a constant uniform magnetic field in $\mathbb{R}^{2n}$\ is
described by the Schrodinger operator 
\begin{equation}
H_{B}=-\frac{1}{4}\sum\limits_{j=1}^{n}\left( \left( \partial
_{x_{j}}+iBy_{j}\right) ^{2}+\left( \partial _{y_{j}}-iBx_{j}\right)
^{2}\right) -\frac{n}{2}  \tag{2.1}
\end{equation}%
acting on $L^{2}(\mathbb{R}^{2n},d\mu )$,\ $d\mu $\ being the Lebesgue
measure on $\mathbb{R}^{2n}$,\ $B>0$\ is the strength of the magnetic
field.\ We identify the Euclidean space $\mathbb{R}^{2n}$\ with $\mathbb{C}%
^{n}$\ in the usual way to use the complex structure. The operator $H_{B}$\
can be represented by the operator 
\begin{equation}\label{eq22}
\Delta _{B}:=e^{\frac{1}{2}B\left\vert z\right\vert ^{2}}H_{B}\text{ }e^{-%
\frac{1}{2}B\left\vert z\right\vert ^{2}}  \tag{2.2}
\end{equation}%
According to \eqref{eq22}, an arbitrary state $\phi $\ of $L^{2}(%
\mathbb{R}^{2n},d\mu )$\ is represented by the function $Q\left[ \phi \right]
$ of $L^{2}(\mathbb{C}^{n},e^{-B\left\vert z\right\vert ^{2}}d\mu )$ defined
by 
\begin{equation}
Q\left[ \phi \right] \left( z\right) :=e^{\frac{1}{2}B\left\vert
z\right\vert ^{2}}\mathsf{\ }\phi \left( z\right) ,\text{ \ \ }z\in \mathbb{C%
}^{n}.  \tag{2.3}
\end{equation}%
The unitary map $Q$\ is called ground state transformation. \ Explicit
expression for the operator \eqref{eq22} is : 
\begin{equation}\label{eq24}
\Delta _{B}=-\sum\limits_{j=1}^{n}\frac{\partial ^{2}}{\partial
z_{j}\partial \overline{z_{j}}}+B\sum\limits_{j=1}^{n}\overline{z_{j}}\frac{%
\partial }{\partial \overline{z_{j}}}.  \tag{2.4}
\end{equation}%
From now on, we take $B=1$\ in \eqref{eq24} and we consider the
magnetic Laplacian $\widetilde{\Delta }:=\Delta _{1}$\ acting in the Hilbert
space $\mathfrak{H}=$\ $L^{2}(\mathbb{C}^{n},e^{-\left\vert z\right\vert
^{2}}d\mu )$ of complex-valued Gaussian square integrable functions. The
spectrum of $\ \widetilde{\Delta }$\ consists of eigenvalues of infinite
multiplicity of the form (Landau levels): $\epsilon _{m}^{1}:=m,$ \ \ $%
m=0,1,2,...$, and admits a spectral decomposition \cite{AM-jmp}. i.e., an
increasing family of projectors $E_{\lambda }$ satisfying 
\begin{equation}
\mathbf{1}_{\mathfrak{H}}=\int\limits_{-\infty }^{+\infty }dE_{\lambda \text{
\ }}\text{\textsf{\ }and \ \ }\widetilde{\Delta }=\int\limits_{-\infty
}^{+\infty }\lambda dE_{\lambda }  \tag{2.5}
\end{equation}%
where $\mathbf{1}_{\mathfrak{H}}$ is the identity operator, and the
projector operator $E_{\lambda }$ are given by,

\begin{equation}
E_{\lambda }\left[ \phi \right] \left( z\right) =\int\limits_{\mathbb{C}^{n}%
\mathsf{\ }}\Phi \left( \lambda ;z,w\right) \phi \left( w\right)
e^{-\left\vert w\right\vert ^{2}}d\mu \left( w\right) \text{ }  \tag{2.6}
\end{equation}%
with the kernel \ 
\begin{equation}\label{eq28}
\Phi _{\lambda }^{\widetilde{\Delta }}\left( z,w\right) =\pi
^{-n}e^{\left\langle z,w\right\rangle }L_{\lfloor \lambda \rfloor }^{\left(
n\right) }(\left\vert z-w\right\vert ^{2})  \tag{2.7}
\end{equation}%
if $\lambda \geq 0$ and $E_{\lambda }\left[ \phi \right] \left( z\right) =0$%
\ when $\lambda <0.$ \ Here $\left[ a\right] $\ denotes the greatest integer
not exceeding $a$\ and 
\begin{equation}
L_{k}^{\left( \alpha \right) }(x)=\sum\limits_{j=o}^{k}\left( -1\right) ^{j}%
\dbinom{k+\alpha }{k-j}\frac{x^{j}}{j!},\text{ \ }\alpha >-1.  \tag{2.8}
\end{equation}%
is the Laguerre polynomial. \ This spectral decomposition leads to the
resolvent kernel of $\ \widetilde{\Delta }$ ($\cite{AM-jmp},$ p.6939$):$%
\begin{equation}\label{eq210}
R(\zeta ;z,w)=-\pi ^{-n}e^{\left\langle z,w\right\rangle }\Gamma (-\zeta
)\Psi \left( -\zeta ,n,\left\vert z-w\right\vert ^{2}\right)   \tag{2.9}
\end{equation}%
where $\Psi $ is the Tricomi function which can be defined ($\cite{Tri},$p.
56) as 
\begin{equation}
\Psi \left( a,c;\xi \right) =\frac{\Gamma \left( 1-c\right) }{\Gamma \left(
a-c+1\right) }._{1}\digamma _{1}\left( a,c;\xi \right) +\frac{\Gamma \left(
c-1\right) }{\Gamma \left( a\right) }\xi ^{1-c}._{1}\digamma _{1}\left(
a,c;\xi \right)   \tag{2.10}
\end{equation}%
where 
\begin{equation}\label{eq212}
\,{_{1}\digamma _{1}}\left( a,c;\xi \right) =\frac{\Gamma (c)}{\Gamma (a)}%
\sum\limits_{j=0}^{+\infty }\frac{\Gamma (a+j)}{\Gamma (c+j)}\frac{\xi ^{j}}{%
j!}  \tag{2.11}
\end{equation}%
is the confluent hypergeometric function defined with the usual condition $%
c\neq 0,-1,-2,\cdots $.\newline

\bigskip 

\textbf{Remark 2.1. }For\textbf{\ }$n=1,$ which corresponds to the planar
Landau Hamiltonian, we should note that, up to a unitary transformation, the
kernel function $\Phi _{\lambda }^{\widetilde{\Delta }}\left( z,w\right) $
in \eqref{eq28} appeared first in  (\cite{DMM}, Eq. (4$),$p.133) as
the Dirac density matrix (with the notation $g\left( \mathbf{r}_{2},\mathbf{r%
}_{1};E\right) $, $\mathbf{r}_{2}\equiv z,\mathbf{r}_{1}\equiv w,E\equiv
\lambda ).$

\medskip 

\textbf{Proposition 2.1}. \textit{The integrated density of states for the
magnetic Laplacian }$\ \widetilde{\Delta }$\textit{\ is given by} 
\begin{equation}
\mathcal{N}_{\widetilde{\Delta }}\left( \lambda \right) =\frac{\left( \lfloor \lambda \rfloor +n\right) !}{\pi ^{n}\lfloor \lambda \rfloor !n!},\text{ }%
n\geq 1.  \tag{2.12}
\end{equation}

The proof is based on the use of \eqref{eq15} and \eqref{eq28} 
\begin{equation}
\mathcal{N}_{\widetilde{\Delta }}\left( \lambda \right) =\frac{1}{Vol\left(
\Omega \right) }\int\limits_{\Omega }\pi ^{-n}e^{\left\langle
z,z\right\rangle }L_{\lfloor \lambda \rfloor }^{\left( n\right)
}(0)e^{-\left\vert z\right\vert ^{2}}d\mu \left( z\right)   \tag{2.13}
\end{equation}%
\begin{equation}\label{eq215}
=\frac{1}{Vol\left( \Omega \right) }\pi ^{-n}L_{\lfloor \lambda \rfloor
}^{\left( n\right) }(0)\int\limits_{\Omega }d\mu \left( x\right) =\pi
^{-n}L_{\lfloor \lambda \rfloor }^{\left( n\right) }(0),  \tag{2.14}
\end{equation}%
together with the fact that 
\begin{equation}
L_{j}^{\left( \alpha \right) }(0)=\frac{\Gamma \left( j+\alpha +1\right) }{%
j!\Gamma \left( \alpha +1\right) }.  \tag{2.15}
\end{equation}

\bigskip \smallskip 

\textbf{Remark 2.2. }When $n=1,$ the expression \eqref{eq212}
reduces to $\mathcal{N}_{\widetilde{\Delta }}\left( \lambda \right) =\pi
^{-1}\left( 1+\lfloor \lambda \rfloor \right) $ which is corresponds (by
changing notations as : $1+2\lambda \equiv E/B$ and $l\equiv \epsilon
_{m}^{1}=m$ ) to well-known IDS for the Landau Hamiltonian on $\mathbb{R}^{2}
$ as pointed out by S. Nakamura $\left( \cite{Nak}\text{, p.149}\right) $
who rederived it by proving that it was independent of the choice of the
boundary condition.

\section{Link with the Heisenberg sublaplacian}

The Heisenberg group $\mathbb{H}^{n}$ is the nilpotent Lie group whose
underlying manifold is $\mathbb{C}^{n}\times \mathbb{R}$ with coordinates $%
\left( z_{1},...,z_{n};\tau \right) =\left( z,\tau \right) $ and whose group
law is 
\begin{equation}
\left( z,\tau \right) .\left( w,s\right) =z+w,\tau +s+2\func{Im}\left\langle
z,w\right\rangle ,  \tag{3.1}
\end{equation}%
where $\left\langle z,w\right\rangle =\sum_{1}^{n}z_{j}\overline{w_{j}}.$
Letting $z=x+iy,$ then, $x_{1},...,x_{n},y_{1},...,y_{n},\tau $ are real
coordinates on $\mathbb{H}^{n}$. We set 
\begin{equation}
\frac{\partial }{\partial z_{j}}=\frac{1}{2}\left( \frac{\partial }{\partial
x_{j}}-i\frac{\partial }{\partial y_{j}}\right) ,\frac{\partial }{\partial 
\overline{z}_{j}}=\frac{1}{2}\left( \frac{\partial }{\partial x_{j}}+i\frac{%
\partial }{\partial y_{j}}\right) .  \tag{3.2}
\end{equation}%
The operator 
\begin{equation}
\Delta _{\mathbb{H}^{n}}:=\sum\limits_{j=0}^{n}\left[ -\frac{\partial ^{2}}{%
\partial z_{j}\partial \overline{z}_{j}}-\left\vert z_{j}\right\vert ^{2}%
\frac{\partial ^{2}}{\partial t^{2}}+i\frac{\partial }{\partial \tau }\left(
z_{j}\frac{\partial }{\partial z_{j}}-\overline{z}_{j}\frac{\partial }{%
\partial \overline{z}_{j}}\right) \right]   \tag{3.3}
\end{equation}%
is left-invariant and subelleptic of order $\frac{1}{2}$ at each $\left(
z,\tau \right) \in \mathbb{H}^{n}$ $(\cite{Foll},$ p.$374)$ \textit{.}

\medskip

Note that operators $\Delta _{\mathbb{H}^{n}}$ and $\widetilde{\Delta }$ \
may be intertwined in the following way by the relation $\left( \cite{AskMou}%
\text{, p}.4\right) :$ 
\begin{equation}\label{eq34}
\frac{1}{2\lambda }\mathcal{\tau }_{\lambda }\circ \left( \frac{1}{2}%
\mathcal{F}\left[ \Delta _{\mathbb{H}^{n}}\right] -n\lambda \right) \circ
\left( \mathcal{\tau }_{\lambda }\right) ^{-1}=\widetilde{\Delta }  \tag{3.4}
\end{equation}%
where $\tau _{\lambda }$ \ is thr map $L^{2}(\mathbb{C}^{n},d\mu
)\rightarrow L^{2}(\mathbb{C}^{n},e^{-\left\vert z\right\vert ^{2}}d\mu )$
defined by 
\begin{equation}
\mathcal{\tau }_{\lambda }\left[ f\right] \left( z\right) :=2^{-n}\lambda
^{-n}e^{\frac{1}{2}\left\vert z\right\vert ^{2}}f\left( \frac{1}{\sqrt{%
2\lambda }}z\right) ,\text{ \ }z\in \mathbb{C}^{n}  \tag{3.5}
\end{equation}%
and $\mathcal{F}:$ $L^{2}(\mathbb{H}_{n}\ ,d\mu \left( z\right)
dt)\rightarrow $ $L^{2}(\mathbb{C}^{n},d\mu )$ denotes the partial Fourier
transform with respect to the variable $t$, defined by 
\begin{equation}
\mathcal{F}\left[ \varphi \right] \left( z,\lambda \right) =\frac{1}{\sqrt{%
\pi }}\int\limits_{\mathbb{R}}\varphi \left( z,t\right) e^{-i\lambda
t}dt=\lim_{\rho \rightarrow +\infty }\int\limits_{\left\vert t\right\vert
\leq \rho }\varphi \left( z,t\right) e^{-i\lambda t}dt  \tag{3.6}
\end{equation}%
which should be understood in the $L^{2}\left( \mathbb{R}\right) $-norm
sense.\ \ \ \ \ \ \ \ \ \ 

\bigskip

The connection given by \eqref{eq34} together with the expression \eqref{eq210} of the resolent kernel of \ the magnetic Laplacian $%
\widetilde{\Delta }$ enables to establish that for $\zeta \in \mathbb{C}$
with $\func{Re}\zeta <0,$the resolvent operator of $\Delta _{\mathbb{H}^{n}}$
has the form ($\cite[p.4, Theorem\, 3.2]{AskMou}):$%
\begin{equation}
\left( \zeta -\Delta _{\mathbb{H}_{n}}\right) ^{-1}\left[ \varphi \right]
\left( z,t\right) =\int\limits_{\mathbb{H}_{n}}\mathcal{R}\left( \zeta
;\left( z,\tau \right) ,\left( w,s\right) \right) \varphi \left( w,s\right)
d\nu \left( w,s\right)   \tag{3.5}
\end{equation}%
where 
\begin{equation}\label{eq36}
\mathcal{R}\left( \zeta ;\left( z,\tau \right) ,\left( w,s\right) \right)
=-2\pi ^{-n-\frac{1}{2}}\int\limits_{0}^{+\infty }x^{n-1}\Gamma \left( \frac{%
n}{2}-\frac{\zeta }{2x}\right) \Psi \left( \frac{n}{2}-\frac{\zeta }{2x}%
,n,2x\left\vert z-w\right\vert ^{2}\right)   \tag{3.6}
\end{equation}%
\begin{equation*}
\times \exp \left( -x\left\vert z-w\right\vert ^{2}\right) \cos \left(
x\left( \tau -s\right) +2x\func{Im}\left\langle z,w\right\rangle \right) dx.
\end{equation*}%
\ Here $d\nu $ denotes the ordinary Lebesgue measure on $\mathbb{R}%
^{2n+1}\equiv \mathbb{H}_{n},$ which is both left and right invariant Haar
measure on the group $\mathbb{H}_{n}$.

\medskip

We should note, in passing, that in the limit case $\zeta =0,$ Eq. \eqref{eq36} reduces to

\begin{equation}\label{eq37}
-\mathcal{R}\left( 0,\left( z,\tau \right) ,\left( w,s\right) \right)
=G_{0}^{F}\left( \left( z,\tau \right) \circ \left( w,s\right) ^{-1}\right) 
\tag{3.7}
\end{equation}%
for any $\left( z,\tau \right) ,\left( w,s\right) \in \mathbb{H}_{n}$, where 
\begin{equation}\label{eq38}
G_{0}^{F}\left( u,t\right) =c_{n}\left\vert \left( u,t\right) \right\vert
^{-2n},\text{ \ \ }\left( u,t\right) \in \mathbb{H}_{n}  \tag{3.8}
\end{equation}%
with $\left\vert \left( u,t\right) \right\vert =\left( \left\vert
u\right\vert ^{4}+t^{2}\right) ^{\frac{1}{4}}$ being the homogeneous norm on 
$\mathbb{H}_{n}.$ 
\begin{equation}\label{eq39}
c_{n}=\left[ n\left( n+1\right) \int\limits_{\mathbb{H}_{n}}\left\vert
u\right\vert ^{2}\left( \left\vert u\right\vert ^{4}+t^{2}+1\right) ^{-\frac{%
1}{2}\left( n+4\right) }d\nu \left( u,t\right) \right] ^{-1}  \tag{3.9}
\end{equation}%
This means that we may recover the fundamental solution for $\Delta _{%
\mathbb{H}_{n}}$ with source at $\left( 0,0\right) $ as 
\begin{equation}
\left( \Delta _{\mathbb{H}_{n}}\left[ \phi \right] ,G_{0}^{F}\right) =\phi
\left( 0\right) ,\text{ }\phi \in \mathcal{C}_{0}^{\infty }\left( \mathbb{H}%
_{n}\right)   \tag{3.10}
\end{equation}%
in agreement with result $(\cite{Foll},$Theorem 2, p. 375$)$ of Folland who
used an analogous fact to $\left\Vert x\right\Vert ^{2-d}$ being (a constant
multiple of) the fundamental solution of the Laplacian on $\mathbb{R}^{d}$ \
with sourece at the origine $0$. In this respect, Appendix A contains
detailed calculations concerning the limit case $\zeta =0$ showing that the
fundamental solution $G_{0}^{F}\left( z,\tau \right) $ in \eqref{eq38} may also be derived from the resolvent kernel  \eqref{eq36} of $%
\Delta _{\mathbb{H}_{n}}\left[ \phi \right] $. In particular, this provides $%
G_{0}^{F}\left( z,\tau \right) $ with the integral 
\begin{equation}
G_{0}^{F}\left( z,\tau \right) =\pi ^{-n-1}2^{n+1}\Gamma \left( \frac{n}{2}%
\right) \int_{0}^{+\infty }x^{n-1}e^{-x|z|^{2}}\Psi \left( \frac{n}{2}%
,n;2x|z|^{2}\right) \cos \left( \tau x\right) dx  \tag{3.11}
\end{equation}%
in terms of the Tricomi $\Psi $-function as mentionned by \cite{GM}$.$

\section{Integrated density of states for $\Delta _{\mathbb{H}_{n}}$}

By using the following summation formula (\cite{Tri}, p. 92): 
\begin{equation}
\Gamma \left( a\right) \Psi \left( a,c;u\right) =\sum\limits_{j=0}^{+\infty }%
\frac{1}{\left( j+a\right) }L_{j}^{\left( c-1\right) }\left( u\right)  
\tag{4.1}
\end{equation}%
in the Abel sense, the authors (\cite{AskMou}, Proposition 5.1) have
obtained the following form for the resolent kernel \eqref{eq36} as 
\begin{equation}
\mathcal{R}\left( \zeta ;\left( z,\tau \right) ,\left( w,s\right) \right)
=\int\limits_{0}^{+\infty }\Phi _{\lambda }^{\Delta _{\mathbb{H}_{n}}}\left(
\left( z,\tau \right) ,\left( w,s\right) \right) \frac{1}{\zeta -\lambda }%
d\lambda   \tag{4.2}
\end{equation}%
where 
\begin{equation}\label{eq43}
\Phi _{\lambda }^{\Delta _{\mathbb{H}_{n}}}\left( \left( z,\tau \right)
,\left( w,s\right) \right) =\lambda ^{n}\pi ^{-n-\frac{1}{2}%
}\sum\limits_{j=0}^{+\infty }\frac{\exp \left( \frac{-\rho \lambda }{2\left(
2j+n\right) }\right) }{\left( 2j+n\right) ^{n+1}}L_{j}^{\left( n-1\right)
}\left( \frac{\lambda \rho }{2j+n}\right) \cos \left( \frac{\tau \lambda }{%
2\left( 2j+n\right) }\right)   \tag{4.3}
\end{equation}%
with $\rho =\left\vert z-w\right\vert ^{2}$ and $\tau =t-s+2\func{Im}%
\left\langle z,w\right\rangle $, is the kernel of the spectral density
operator $\frac{dE_{\lambda }}{d\lambda }$ as proved in \cite{AskMou} where
the distributional sense is discussed.

\bigskip

\bigskip \textbf{Proposition\medskip\ 4.1.} \textit{The integtated density
of states for }$\Delta _{\mathbb{H}_{n}}$\textit{\ is given by }%
\begin{equation}
\mathcal{N}\left( \lambda \right) =\gamma _{n}\lambda ^{n}\mathit{\ } 
\tag{4.4}
\end{equation}%
\textit{where }$\gamma _{n}>0$\textit{\ is the constant given by}%
\begin{equation}
\gamma _{n}=\frac{1}{\pi ^{n+\frac{1}{2}}\Gamma \left( n\right) }%
\sum\limits_{j=0}^{+\infty }\frac{\Gamma \left( j+n\right) }{j!\left(
2j+n\right) ^{n+1}}  \tag{4.5}
\end{equation}

\bigskip \medskip

\textbf{Proof.} \textbf{\ }Let $\Omega \subset \mathbb{H}_{n}$\ be a domain
of finite volume. According to \eqref{eq15} , the IDS can be
expressed as 
\begin{equation}\label{eq46}
\mathcal{N}\left( \lambda \right) :=\frac{1}{Vol\left( \Omega \right) }%
\int\limits_{\Omega }\Phi _{\lambda }^{\Delta _{\mathbb{H}_{n}}}\left(
\left( z,\tau \right) ,\left( w,s\right) \right) d\mu \left( w,s\right) . 
\tag{4.6}
\end{equation}%
By substiting $\Phi _{\lambda }^{\Delta _{\mathbb{H}_{n}}}\left( \left(
z,\tau \right) ,\left( w,s\right) \right) $ by its expression \eqref{eq43} into \eqref{eq46}, Eq. \eqref{eq46} reduces to 
\begin{equation}
\mathcal{N}\left( \lambda \right) =\frac{\lambda ^{n}}{\pi ^{n+\frac{1}{2}}}%
\sum\limits_{j=0}^{+\infty }\frac{1}{\left( 2j+n\right) ^{n+1}}L_{j}^{\left(
n-1\right) }\left( 0\right)  \tag{4.7}
\end{equation}%
Finally, we use \eqref{eq215} in order to get the exact value of $%
L_{j}^{\left( n-1\right) }\left( 0\right) $.
\newpage
\bigskip

\begin{center}
\textbf{Appendix A}

\medskip 
\end{center}

To prove \eqref{eq37}, we first recall \ \cite[p.4, Theorem 3.2]%
{AskMou}, the resolvent operator of $\Delta _{\mathbb{H}%
_{n}}$ has the form 
\begin{equation}\label{A1}
\left( \zeta -\Delta _{\mathbb{H}_{n}}\right) ^{-1}\left[ \varphi \right]
\left( z,t\right) =\int_{\mathbb{H}^{n}}\mathcal{R}\left( \zeta ;\left(
z,\tau \right) ,\left( w,s\right) \right) \varphi \left( w,s\right) d\mu
(w,s),\text{ }\func{Re}\zeta <0,  \tag{A1}
\end{equation}%
where 
\begin{equation}
\mathcal{R}\left( \zeta ;\left( z,\tau \right) ,\left( w,s\right) \right) =%
\frac{-2^{n}}{\pi ^{n+\frac{1}{2}}}\int_{0}^{+\infty }x^{n-1}\Gamma \left( 
\frac{n}{2}-\frac{\zeta }{2x}\right) \Psi \left( \frac{n}{2}-\frac{\zeta }{2x%
},n;2x\left\vert z-w\right\vert ^{2}\right)   \tag{A2}
\end{equation}%
\begin{equation*}
\times \exp \left( -x\left\vert z-w\right\vert ^{2}\right) \cos \left(
x\left( \tau -s\right) +2x\func{Im}\left\langle z,w\right\rangle \right) dx.
\end{equation*}%
\newline
At the limit $\zeta \rightarrow 0$ in \eqref{A1}, we obtain a right
inverse of $\Delta _{\mathbb{H}_{n}}$ as 
\begin{equation}
\Delta _{\mathbb{H}_{n}}^{-1}\left[ \varphi \right] \left( z,\tau \right)
=\int_{\mathbb{H}^{n}}\mathcal{R}\left( 0;\left( z,\tau \right) ,\left(
w,s\right) \right) \varphi \left( w,s\right) d\nu (w,s),  \tag{A3}
\end{equation}%
or equivalently 
\begin{equation}
\mathcal{R}_{0}:=-\mathcal{R}\left( 0;(z,\tau ),(w,s)\right)   \tag{A4}
\end{equation}%
is a Green kernel for $\Delta _{\mathbb{H}_{n}}$ as pointed out in \cite[%
p.6, Remark 3.3]{AskMou}. Now, to establish a connection between the
integral kernel $\mathcal{R}_{0}$ and the Folland's fundamental solution, we
proceed by computing the integral 
\begin{equation}\label{A5}
\mathcal{R}_{0}=2^{n}\pi ^{-n-\frac{1}{2}}\int_{0}^{+\infty }x^{n-1}\Gamma
\left( \frac{n}{2}\right) \Psi \left( \frac{n}{2},n;2x|z-w|^{2}\right)  
\tag{A5}
\end{equation}%
\begin{equation*}
\times \exp \left( -x\left\vert z-w\right\vert ^{2}\right) \cos \left(
x\left( \tau -s\right) +2x\func{Im}\left\langle z,w\right\rangle \right) dx.
\end{equation*}%
For this, let us rewrite \eqref{A5} as 
\begin{equation}\label{A6}
\mathcal{R}_{0}=2^{n}\pi ^{-n-\frac{1}{2}}\int_{0}^{+\infty }x^{n-1}\left[
\Gamma \left( \frac{n}{2}\right) \Psi \left( \frac{n}{2},n;\mu x\right) %
\right] e^{-\frac{\mu }{2}x}\cos \left( x\theta \right) dx,  \tag{A6}
\end{equation}%
where we have set $\mu :=2|z-w|^{2}$ and $\theta =\left( \tau -s\right) +2%
\func{Im}\left\langle z,w\right\rangle $. By using the integral
representation of $\Psi $-function (\cite{Tri}, p. 92) : 
\begin{equation}
\Gamma \left( a\right) \Psi \left( a,c;u\right) =\int_{0}^{+\infty
}e^{-at}\exp \left( -\frac{u}{e^{t}-1}\right) \left( 1-e^{-t}\right) ^{-c}dt
\tag{A7}
\end{equation}%
for parameters $a=\frac{n}{2}$, $c=n$ and $u=\mu x$, the integral \eqref{A6} reads 
\begin{equation}\label{A8}
\mathcal{R}_{0}=\frac{2^{n}}{\pi ^{n+\frac{1}{2}}}\int_{0}^{+\infty }x^{n-1}%
\left[ \int_{0}^{+\infty }e^{-\frac{n}{2}t}\exp \left( -\frac{\mu x}{e^{t}-1}%
\right) \left( 1-e^{-t}\right) ^{-n}dt\right] e^{-\frac{\mu }{2}x}\cos
\left( x\theta \right) dx.  \tag{A8}
\end{equation}%
Now, applying the identity (\cite{GR}, p.498): 
\begin{equation}
\int_{0}^{+\infty }x^{\nu -1}e^{-\alpha x}\cos (\theta x)dx=\Gamma \left(
\nu \right) \left( \alpha ^{2}+\theta ^{2}\right) ^{-\frac{\nu }{2}}\cos
\left( \nu \arctan \left( \frac{\theta }{\alpha }\right) \right) ,\func{Re}%
\nu >0,\func{Re}\alpha >|\func{Im}(\theta )|  \tag{A9}
\end{equation}%
for parameters $\nu =n$ and $\alpha =\frac{\mu }{2}\coth \left( \frac{t}{2}%
\right) $ and Eq. \eqref{A8} reads 
\begin{equation}
\mathcal{R}_{0}=2^{n}\Gamma (n)\pi ^{-n-\frac{1}{2}}\int_{0}^{+\infty }\frac{%
e^{-\frac{n}{2}t}}{\left( 1-e^{-t}\right) ^{n}}\left( \left( \frac{\mu }{2}%
\coth \left( \frac{t}{2}\right) \right) ^{2}+\theta ^{2}\right) ^{-\frac{n}{2%
}}\cos \left( n\arctan \left( \frac{\theta }{\frac{\mu }{2}\coth \left( 
\frac{t}{2}\right) }\right) \right) dt  \tag{A10}
\end{equation}%
\begin{equation}\label{A11}
=\frac{\Gamma (n)}{\pi ^{n+\frac{1}{2}}}\int_{0}^{+\infty }\left( \sinh
\left( \frac{t}{2}\right) \right) ^{-n}\left( \frac{\mu }{2}\coth \left( 
\frac{t}{2}\right) \right) ^{-n}\left( 1+\left( \frac{\theta }{\frac{\mu }{2}%
\coth \left( \frac{t}{2}\right) }\right) ^{2}\right) ^{-\frac{n}{2}}\cos
\left( n\arctan \left( \frac{\theta }{\frac{\mu }{2}\coth \left( \frac{t}{2}%
\right) }\right) \right) dt  \tag{A11}
\end{equation}%
Making the change variable $\rho =\frac{\theta }{\frac{\mu }{2}\coth \frac{t%
}{2}},$ \eqref{A11} becomes 
\begin{equation}\label{A12}
\mathcal{R}_{0}=\frac{2^{n}\Gamma (n)}{\mu ^{n-1}\theta \pi ^{n+\frac{1}{2}}}%
\int_{0}^{\frac{2\theta }{\mu }}\left( 1-\left( \frac{\mu }{2\theta }\right)
^{2}\rho ^{2}\right) ^{\frac{n}{2}-1}\left( 1+\rho ^{2}\right) ^{-\frac{n}{2}%
}\cos \left( n\arctan (\rho )\right) d\rho .  \tag{A12}
\end{equation}%
A second change of variable, $\arctan \left( \rho \right) =u$, shows that \eqref{A12} can be reduced to 
\begin{equation}
\mathcal{R}_{0}=\frac{2^{n}\Gamma (n)}{\mu ^{n-1}\theta \pi ^{n+\frac{1}{2}}}%
\int_{0}^{\arctan (\frac{2\theta }{\mu })}\left( \cos ^{2}(u)-\left( \frac{%
\mu }{2\theta }\right) ^{2}\sin ^{2}u\right) ^{\frac{n}{2}-1}\cos (nu)du. 
\tag{A13}
\end{equation}%
\begin{equation}\label{A14}
=\frac{2^{\frac{n}{2}+1}\Gamma (n)}{\mu ^{n-1}\theta \pi ^{n+\frac{1}{2}}}%
\left( 1+\left( \frac{\mu }{2\theta }\right) ^{2}\right) ^{\frac{n}{2}%
-1}\int_{0}^{\arctan (\frac{2\theta }{\mu })}\left( \cos 2u-\frac{1-\left( 
\frac{2\theta }{\mu }\right) ^{2}}{1+\left( \frac{2\theta }{\mu }\right) ^{2}%
}\right) ^{\frac{n}{2}-1}\cos (nu)du.  \tag{A14}
\end{equation}%
Setting $\beta =\left( \frac{2\theta }{\mu }\right) ^{2}$ and making the
change of variable $2u=\kappa $ in \eqref{A14} gives that 
\begin{equation}\label{A15}
\mathcal{R}_{0}=\frac{2^{\frac{n}{2}}\Gamma (n)}{\mu ^{n-1}\theta \pi ^{n+%
\frac{1}{2}}}\left( \frac{\beta +1}{\beta }\right) ^{\frac{n}{2}%
-1}\int_{0}^{2\arctan (\frac{2\theta }{\mu })}\left( \cos \kappa -\frac{%
1-\beta }{1+\beta }\right) ^{\frac{n}{2}-1}\cos \left( \frac{n}{2}\kappa
\right) d\kappa .  \tag{A15}
\end{equation}%
Now, by putting 
\begin{equation}\label{A16}
\varepsilon =\arccos \left( \frac{1-\beta }{1+\beta }\right) =2\arctan
\left( \frac{2\theta }{\mu }\right) .  \tag{A16}
\end{equation}%
We are now in position to apply the identity (\cite{GR}, p.406): 
\begin{equation}\label{A17}
\int_{0}^{\varepsilon }\left( \cos (x)-\cos (\varepsilon )\right) ^{\nu -%
\frac{1}{2}}\cos (ax)dx=\sqrt{\frac{\pi }{2}}\left( \sin (\varepsilon
)\right) ^{\nu }\Gamma \left( \nu +\frac{1}{2}\right) P_{a-\frac{1}{2}%
}^{-\nu }\left( \cos (\varepsilon )\right) ,\text{ }\func{Re}\nu >-\frac{1}{2%
},a>0,\text{ }0<\varepsilon <\pi   \tag{A17}
\end{equation}%
where 
\begin{equation}
P_{\lambda }^{\nu }(x):=\frac{1}{\Gamma (1-\nu )}\left( \frac{x+1}{x-1}%
\right) ^{\nu /2}{_{2}\digamma _{1}}\left( -\lambda ,\lambda +1;1-\nu ;\frac{%
1-x}{2}\right)   \tag{A18}
\end{equation}%
denotes the Legendre function of the first kind (\cite{GR}, p. 959)
In our setting, the integral in \eqref{A17} reads 
\begin{equation}\label{A19}
\int_{0}^{2\arctan (\frac{2\theta }{\mu })}\left( \cos (\kappa )-\frac{%
1-\beta }{1+\beta }\right) ^{\frac{n-1}{2}-\frac{1}{2}}\cos \left( \frac{n}{2%
}\kappa \right) d\kappa =\sqrt{\frac{\pi }{2}}\left( \sin (\varepsilon
)\right) ^{\frac{n-1}{2}}\Gamma \left( \frac{n}{2}\right) P_{\frac{n-1}{2}%
}^{-\left( \frac{n-1}{2}\right) }\left( \cos \varepsilon \right) .  \tag{A19}
\end{equation}%
Hence, using the Gegenbauer representation (\cite{GR}, p.959): 
\begin{equation}
P_{\sigma }^{-\sigma }\left( \cos (\varepsilon )\right) =\frac{1}{\Gamma
\left( 1+\sigma \right) }\left( \frac{1}{2}\sin (\varepsilon )\right)
^{\sigma }  \tag{A20}
\end{equation}%
with $\sigma =\frac{n-1}{2}$, we can write the right hand side in \eqref{A19} as 
\begin{equation}
\sqrt{\frac{\pi }{2}}\left( \sin (\varepsilon )\right) ^{\frac{n-1}{2}%
}\Gamma \left( \frac{n}{2}\right) P_{\frac{n-1}{2}}^{-\left( \frac{n-1}{2}%
\right) }\left( \cos (\varepsilon )\right) =\left( \frac{1}{2}\right) ^{%
\frac{n}{2}}\sqrt{\pi }\frac{\Gamma \left( \frac{n}{2}\right) }{\Gamma
\left( \frac{n+1}{2}\right) }\left( 1-\cos ^{2}(\varepsilon )\right) ^{\frac{%
n-1}{2}}.  \tag{A21}
\end{equation}%
Returning back to \eqref{A15}, keeping in mind the expression of $%
\varepsilon $ given through \eqref{A16}, we obtain that 
\begin{equation}
\mathcal{R}_{0}=\frac{2^{n-1}}{\mu ^{n-1}\theta \pi ^{n}}\frac{\Gamma
(n)\Gamma \left( \frac{n}{2}\right) }{\Gamma \left( \frac{n+1}{2}\right) }%
\left( \frac{\beta +1}{\beta }\right) ^{-\frac{n}{2}}.  \tag{A22}
\end{equation}%
Now, replacing $\beta $ by $\left( \frac{2\theta }{\mu }\right) ^{2}$ and
using the expressions of $\mu $ and $\theta $, we arrive at 
\begin{equation}
\mathcal{R}_{0}=\frac{\Gamma (n)\Gamma \left( \frac{n}{2}\right) }{\pi
^{n}\Gamma \left( \frac{n+1}{2}\right) }\left( \left( |z-w|^{2}\right)
^{2}+\left( \left( \tau -s\right) +2\Im \left\langle z,w\right\rangle
\right) ^{2}\right) ^{-\frac{n}{2}}  \tag{A23}
\end{equation}%
\begin{equation}
=\frac{2^{n-1}\Gamma ^{2}\left( \frac{n}{2}\right) }{\pi ^{n+\frac{1}{2}}}%
\left( \left( |z-w|^{2}\right) ^{2}+\left( \left( \tau -s\right) +2\Im
\left\langle z,w\right\rangle \right) ^{2}\right) ^{-\frac{n}{2}}.  \tag{A24}
\end{equation}%
The last equality follows using Legendre's duplication formula (\cite{GR}, p.896): 
\begin{equation}
\Gamma \left( \xi \right) \Gamma \left( \xi +\frac{1}{2}\right) =2^{1-2\xi }%
\sqrt{\pi }\Gamma (2\xi )  \tag{25}
\end{equation}%
for $\xi =\frac{n}{2}$. Therefore, we assert that 
\begin{equation}\label{A26}
\mathcal{R}_{0}=\frac{2^{n-1}\Gamma ^{2}\left( \frac{n}{2}\right) }{\pi ^{n+%
\frac{1}{2}}c_{n}}G_{0}^{F}\left( \left( z,\tau \right) \circ \left(
w,s\right) ^{-1}\right) ,  \tag{A26}
\end{equation}%
where the constant $c_{n}$ is as in \eqref{eq39}. In particular, for 
$\left( w,s\right) =\left( 0,0\right) $, Eq. \eqref{A26} reduces
further to 
\begin{equation}
\mathcal{R}_{0}=\frac{2^{n-1}\Gamma ^{2}\left( \frac{n}{2}\right) }{\pi ^{n+%
\frac{1}{2}}c_{n}}G_{0}^{F}\left( z,\tau \right) =\frac{\sqrt{\pi }}{2}%
G_{0}^{F}\left( z,\tau \right) .  \tag{A27}
\end{equation}%
This completes the proof of \eqref{eq37}.

\bigskip \medskip

\begin{quote}
\bigskip 
\end{quote}

\end{document}